\documentclass[12pt,a4paper,reqno]{amsart}
\usepackage{amssymb,graphpap}
\usepackage{amscd}
\theoremstyle{definition}

\def\rb{\rbrack}
\def\lb{\lbrack}
\def\la{\langle}
\def\ra{\rangle}

\def\cA{\mathcal{A}}
\def\cF{\mathcal{F}}
\def\cP{\mathcal{P}}
\def\cB{\mathcal{B}}

\def\HH{\mathcal{H}}
\def\CC{{\Bbb C}}
\def\NN{{\Bbb N}}

\def\ZZ{{\Bbb Z}}

\def\ff{\varphi}
\def\folgt{\qquad\Longrightarrow\qquad}
   
\def\folgtohne{\Longrightarrow}

\def\entspricht{\qquad\hat =\qquad}
\def\entsprichtohne{\hat =}
\def\cupdot{\dot\cup}

\def\tr{\text{tr}}

\def\qundq{\qquad\text{and}\qquad}
\def\FPR{\Theta}
\def\zz{\overset {!}{=}}

\def\distr{\text{distr}}
\def\eins{{\mathbf 1}}
\def\null{{\mathbf 0}}

\begin{document}  

\title{Free calculus}
\author{Roland Speicher}
\address{\noindent Department of Mathematics and Statistics,
Queen's University, Kingston, ON K7L 3N6, Canada}
\email{speicher@mast.queensu.ca}

\maketitle

\section{{\bf Basic definitions and facts}}

\subsection*{}
Let me first recall the basic definitions 
and some fundamental realizations of 
freeness.
For a more extensive review, I refer to the course of Biane or to
the books \cite{VDN,V6} (see also the survey \cite{V5}).
 
\subsection{Definitions}
1)  A {\bf (non-commutative) probability space}
consists of a pair $(\cA,\ff)$, where
\begin{itemize}
\item
$\cA$ is a unital algebra
\item
$\ff:\cA\to\CC$ is a unital linear functional, i.e. in particular
$\ff(1)=1$
\end{itemize}
2) Unital subalgebras $\cA_1,\dots,\cA_n\subset\cA$ are called
{\bf free}, if we have
$$\ff(a_1\dots a_k)=0$$
whenever
\begin{gather*}
a_i\in\cA_{j(i)}\qquad (i=1,\dots,k)\\
j(1)\not=j(2)\not=\dots\not=j(k)\\
\ff(a_i)=0\qquad (i=1,\dots,k)
\end{gather*}
3) Random variables $x_1,\dots,x_n\in\cA$ are
called {\bf free}, if $\cA_1,\dots,\cA_n$ are free, where
$\cA_i$ is the unital algebra generated by $x_i$.

\subsection*{}
A canonical realization of free random variables is given
on the full Fock space.

\subsection{Definitions}
Let $\HH$ be a Hilbert space.
\\
1) The {\bf full Fock space over $\HH$} is the Hilbert space
$$\cF(\HH):=\CC\Omega\oplus\bigoplus_{n\geq 1}\HH^{\otimes n},$$
where $\Omega$ is a distinguished unit vector, called {\bf vacuum}.
\\
2) The {\bf vacuum expectation} is the state
$$A\mapsto \la \Omega,A\Omega\ra.$$
\\
3) For each $f\in\HH$ we define the {\bf 
(left) annihilation
operator $l(f)$} and the {\bf (left) creation operator $l^*(f)$}
by
\begin{align*}
l(f)\Omega&=0\\
l(f)f_1\otimes \dots\otimes f_n&=\la f,f_1\ra f_2\otimes 
\dots\otimes f_n
\end{align*}
and
$$l^*(f)f_1\otimes \dots\otimes f_n=f\otimes f_1\otimes\dots\otimes f_n.$$

One should note that in regard of whether we denote by $l(f)$ the
annihilation or the creation operator we follow the opposite
tradition as Voiculescu.

\subsection{Proposition}
Let $\HH_1$ and $\HH_2$ be Hilbert spaces and put
$\HH:=\HH_1\oplus \HH_2$.
Consider the full Fock space over $\HH$ and the corresponding
creation and annihilation operators $l(f)$ and $l^*(f)$ for
$f\in\HH$.
Put
$$\cA_1:=C^*(l(f)\mid f\in\HH_1),\qquad
\cA_2:=C^*(l(f)\mid f\in\HH_2).$$
Then $\cA_1$ and $\cA_2$ are free with respect to the vacuum expectation.

\subsection*{}
If we replace the $C^*$-algebras by von Neumann algebras or
by $*$-algebras, then the analogue statements are also true. The proof
consists of direct checking for the case of $*$-algebras and
then extending the assertion to uniform or weak closure by approximation
arguments.

\subsection*{}
Whereas freeness is just modelled according to the situation on
the full Fock space -- hence its appearing in this context is not
very surprising -- there is another realization of freeness in a totally
different context: freeness can also be thought of as the mathematical
structure of $N\times N$ random matrices in the limit $N\to\infty$.
We will not need this connection in our considerations, but it is
always good to keep in mind that all our constructions 
also have some meaning
in terms of random matrices.

Such random matrices $X_N$ are $N\times N$ matrices whose entries
are classical random variables and usually one is interested in
the averaged eigenvalue distribution $\text{distr}(X_N)$ of
these matrices corresponding to a state given by the averaged
normalized trace.

The following theorem is due to Voiculescu \cite{V4} (see also
\cite{SpRM}).

\subsection{Theorem}
{\bf 1) Gaussian random matrices}
\\
Let 
$$X_N= (a^{(N)}_{ij})_{i,j=1}^N,\qquad
Y_N= (b^{(N)}_{ij})_{i,j=1}^N$$
be symmetric random $N\times N$ matrices with
\begin{itemize}
\item
$a^{(N)}_{ij}$ ($1\leq i\leq j\leq N$) are mutually independent and normally
distributed with mean zero and variance $1/N$.
\item
$b^{(N)}_{ij}$ ($1\leq i\leq j\leq N$) are mutually independent and normally
distributed with mean zero and variance $1/N$.
\item
all $a^{(N)}_{ij}$ are independent from all $b^{(N)}_{kl}$.
\end{itemize}
Then $X_N$ and $Y_N$ converge in distribution
to a semicircular family 
in the limit $N\to\infty$ with respect to
$$\ff(\cdot):=\langle\frac 1N\tr(\cdot)\rangle_{ensemble}.$$
{\bf 2) randomly rotated matrices}
\\
Let $A_N$ and $B_N$ be symmetric deterministic 
(e.g. diagonal) $N\times N$ matrices with
$$\lim_{N\to\infty}\distr(A_N)=\mu,\qquad
\lim_{N\to\infty}\distr(B_N)=\nu$$
for some compactly supported probability measures $\mu$ and $\nu$.
Let $U_N$ be a random unitary matrix from the ensemble
$$\Omega_N=(\mathcal{U}(N),Haar\  measure).$$
Consider now
$$X_N:=A_N,\qquad Y_N:=U_NB_NU_N^*.$$
Then $X_N$ and $Y_N$ become free 
in the limit $N\to\infty$ with respect to
$$\ff(\cdot):=\langle \frac 1N\tr(\cdot)\rangle_{\Omega_N}.$$

\section{{\bf Combinatorial aspects of
freeness: the concept of cumulants}}

\subsection*{}
`Freeness' of random variables is defined in terms of mixed moments;
namely the defining property is that very special moments (alternating
and centered ones) have to vanish. This requirement is not easy to
handle in concrete calculations. Thus we will present here another
approach to freeness, more combinatorial in nature, which puts the
main emphasis on so called `free cumulants'. These are some
polynomials in the moments which behave much better with respect
to freeness than the moments. The nomenclature comes from classical
probability theory where corresponding objects are also well
known and are usually called `cumulants' or `semi-invariants'.
There exists a combinatorial description of these classical cumulants,
which depends on partitions of sets. In the same way, free cumulants
can also be described combinatorially, the only difference to
the classical case is that one has to replace all partitions
by so called `non-crossing partitions'.

In the case of one random variable,
we will also indicate the relation of this combinatorial description
with the analytical one presented in the course of Biane; namely
our cumulants are in this case
just the coefficients of the $R$-transform of
Voiculescu (in the classical case the cumulants are the coefficients
of the logarithm of the Fourier transform). Thus we will obtain
purely combinatorial proofs of the main results on the $R$-transform.

This combinatorial description of freeness is due to me
\cite{Sp1,Sp2,Sp3} (see also \cite{N}); in a series of joint
papers with A. Nica \cite{NS2,NS3,NS4}
it was pursued very far and yielded a
lot of new results in free probability theory. I will restrict
here mainly to the basic facts; for applications one should
consult the original papers or my survey \cite{Sp4}.
A recent fundamental link between freeness and the representation
theory of the permutation groups $S_n$ in the limit $n\to\infty$,
which rests also on the combinatorial description of freeness, is
due to Biane \cite{B}.

\subsection{Definitions}
A {\bf partition} of the set $S:=\{1,\dots,n\}$ is a
decomposition
$$\pi=\{V_1,\dots,V_r\}$$
of $S$ into disjoint and non-empty sets $V_i$, i.e.
$$V_i\not=\emptyset\quad(i=1,\dots,r)
\qundq S=\cupdot_{i=1}^rV_i.$$
We denote the set of all partitions of $S$ with
$\cP(S)$.
\\
We call the $V_i$ the {\bf blocks} of $\pi$.
\\
For $1\leq p,q\leq n$ we write
$$p\sim_\pi q\qquad\text{if $p$ and $q$ belong to the
same block of $\pi$.}$$
A partition $\pi$ is called {\bf non-crossing}
if the following does not occur: There exist
$1\leq p_1<q_1<p_2<q_2\leq n$ with
$$p_1\sim_\pi p_2\not\sim_\pi q_1\sim_\pi q_2.$$
The set of all non-crossing partitions of
$\{1,\dots,n\}$
is denoted by {\bf NC(n)}.
\\
We denote the `biggest' and the
`smallest' element in $NC(n)$ by $\eins_n$ and $\null_n$,
respectively:
\begin{align*}
\eins_n:&=\{(1,\dots,n)\},\qquad
\null_n:=\{(1),\dots,(n)\}.
\end{align*}

Non-crossing partitions were introduced by Kreweras \cite{Kr}
in a purely combinatorial context without any reference to probability
theory.

\subsection{Examples}
We will also use a graphical notation for our
partitions; the term `non-crossing' will become evident in
such a notation.
\setlength{\unitlength}{0.3cm}
Let
$$S=\{1,2,3,4,5\}.$$
Then the partition
$$\pi=\{(1,3,5),(2),(4)\}
\entspricht
\begin{picture}(4,2)\thicklines
\put(0,0){\line(0,1){2}}
\put(0,0){\line(1,0){4}}
\put(2,0){\line(0,1){2}}
\put(4,0){\line(0,1){2}}
\put(1,1){\line(0,1){1}}
\put(3,1){\line(0,1){1}}
\put(-0.3,2.3){1}
\put(0.7,2.2){2}
\put(1.7,2.2){3}
\put(2.7,2.2){4}
\put(3.7,2.2){5}
\end{picture}
$$
is non-crossing, whereas
$$\pi=\{(1,3,5),(2,4)\}
\entspricht
\begin{picture}(4,2)\thicklines
\put(0,0){\line(0,1){2}}
\put(0,0){\line(1,0){4}}
\put(2,0){\line(0,1){2}}
\put(4,0){\line(0,1){2}}
\put(1,1){\line(0,1){1}}
\put(3,1){\line(0,1){1}}
\put(1,1){\line(1,0){2}}
\put(-0.3,2.3){1}
\put(0.7,2.2){2}
\put(1.7,2.2){3}
\put(2.7,2.2){4}
\put(3.7,2.2){5}
\end{picture}
$$
is crossing.

\subsection{Remarks}
1) In an analogous way, non-crossing partitions $NC(S)$ can be
defined for any linearly ordered set $S$;
of course, we have
$$NC(S_1)\cong NC(S_2)\qquad\text{if}\qquad \#S_1=\#S_2.$$
2) In most cases the following recursive description of
non-crossing partitions is of great use:
a partition $\pi$ ist non-crossing if and only if
at least one block $V\in\pi$ is an interval and
$\pi\backslash V$ is non-crossing; i.e. 
$V\in\pi$ has the form
$$V=(k,k+1,\dots,k+p)\qquad\text{for some $1\leq k\leq n$
and $p\geq 0$, $k+p\leq n$}$$
and we have
$$\pi\backslash V\in NC(1,\dots,k-1,k+p+1,\dots,n)\cong
NC(n-(p+1)).$$
Example: The partition
$$\{(1,10),(2,5,9),(3,4),(6),(7,8)\}
\entspricht
\setlength{\unitlength}{0.3cm}
\begin{picture}(9,3)\thicklines
\put(0,0){\line(0,1){3}}
\put(0,0){\line(1,0){9}}
\put(9,0){\line(0,1){3}}
\put(1,1){\line(0,1){2}}
\put(1,1){\line(1,0){7}}
\put(4,1){\line(0,1){2}}
\put(8,1){\line(0,1){2}}
\put(2,2){\line(0,1){1}}
\put(2,2){\line(1,0){1}}
\put(3,2){\line(0,1){1}}
\put(5,2){\line(0,1){1}}
\put(6,2){\line(0,1){1}}
\put(6,2){\line(1,0){1}}
\put(7,2){\line(0,1){1}}
\put(-0.3,3.3){1}
\put(0.7,3.3){2}
\put(1.7,3.3){3}
\put(2.7,3.3){4}
\put(3.7,3.3){5}
\put(4.7,3.3){6}
\put(5.7,3.3){7}
\put(6.7,3.3){8}
\put(7.7,3.3){9}
\put(8.7,3.3){10}
\end{picture}
$$
can, by successive removal of intervals, be reduced to
$$\{(1,10),(2,5,9)\}\entsprichtohne \{(1,5),(2,3,4)\}$$
and finally to
$$\{(1,5)\}\entsprichtohne \{(1,2)\}.$$
3) By writing a partition $\pi$ in the form
$\pi=\{V_1,\dots,V_r\}$ we will always assume that the
elements within each block $V_i$
are ordered in increasing order.

\subsection{Definition}
Let $(\cA,\ff)$ be a probability space,
i.e. $\cA$ is a unital algebra and
$\ff:\cA\to\CC$ is a unital linear functional.
We define the {\bf (free or non-crossing) cumulants} 
$$k_n:\cA^n\to\CC\qquad (n\in\NN)$$
(indirectly) by the following system of equations:
$$\ff(a_{1}\dots a_{n})=\sum_{\pi\in NC(n)}
k_{\pi}\lb a_{1},\dots,a_{n}\rb\qquad (a_1,\dots,a_n\in\cA),$$
where $k_\pi$ denotes a product of cumulants according to
the block structure of $\pi$:
\begin{multline*}
k_{\pi}\lb a_{1},\dots,a_{n}\rb:=
k_{V_1}\lb a_{1},\dots,a_{n}\rb\dots
k_{V_r}\lb a_{1},\dots,a_{n}\rb\qquad
\text{for $\pi=\{V_1,\dots,V_r\}\in NC(n)$}
\end{multline*}
and
$$
k_{V}\lb a_{1},\dots,a_{n}\rb:=k_{\#V}(a_{v_1},\dots,
a_{v_l})\qquad\text{for}\qquad V=(v_1,\dots,v_l).$$

\subsection{Remarks and Examples}
1) Note: the above equations have the form
$$\ff(a_{1}\dots a_{n})=k_n(a_{1},\dots, a_{n})
+\sum_{\pi\in NC(n)\atop \pi\not=\eins_n} k_\pi
\lb a_1,\dots,a_n\rb,$$
and thus they can be resolved for the
$k_n(a_{1},\dots, a_{n})$
in a unique way.
\\
2) Examples:
\setlength{\unitlength}{0.2cm}
\begin{itemize}
\item
$n=1$
$$\ff(a_1)=k_{\,
\begin{picture}(0,0)\thicklines
\put(0,0){\line(0,1){1}}
\end{picture}\,
}
\lb a_1\rb=k_1(a_1),$$
thus
$$k_1(a_1)=\ff(a_1).$$
\item
$n=2$
\begin{align*}
\ff(a_1a_2)&=k_{\,
\begin{picture}(1,1)\thicklines
\put(0,0){\line(0,1){1}}
\put(0,0){\line(1,0){1}}
\put(1,0){\line(0,1){1}}
\end{picture}\,
}
\lb a_1,a_2\rb +k_{\,
\begin{picture}(1,1)\thicklines
\put(0,0){\line(0,1){1}}
\put(1,0){\line(0,1){1}}
\end{picture}\,
}
\lb a_1,a_2\rb\\
&=k_2(a_1,a_2)+k_1(a_1)k_1(a_2),
\end{align*}
thus
$$k_2(a_1,a_2)=\ff(a_1a_2)-\ff(a_1)\ff(a_2).$$
\item
$n=3$
\begin{align*}
\ff(a_1a_2a_3)&=
k_{\,
\begin{picture}(2,1)\thicklines
\put(0,0){\line(0,1){1}}
\put(0,0){\line(1,0){2}}
\put(1,0){\line(0,1){1}}
\put(2,0){\line(0,1){1}}
\end{picture}\,
}
\lb a_1,a_2,a_3\rb
+k_{\,
\begin{picture}(2,1)\thicklines
\put(0,0){\line(0,1){1}}
\put(1,0){\line(1,0){1}}
\put(1,0){\line(0,1){1}}
\put(2,0){\line(0,1){1}}
\end{picture}\,
}
\lb a_1,a_2,a_3\rb
+k_{\,
\begin{picture}(2,1)\thicklines
\put(0,0){\line(0,1){1}}
\put(0,0){\line(1,0){1}}
\put(1,0){\line(0,1){1}}
\put(2,0){\line(0,1){1}}
\end{picture}\,
}
\lb a_1,a_2,a_3\rb\\
&\quad
+k_{\,
\begin{picture}(2,1)\thicklines
\put(0,0){\line(0,1){1}}
\put(0,0){\line(1,0){2}}
\put(1,0.5){\line(0,1){0.5}}
\put(2,0){\line(0,1){1}}
\end{picture}\,
}
\lb a_1,a_2,a_3\rb
+k_{\,
\begin{picture}(2,1)\thicklines
\put(0,0){\line(0,1){1}}
\put(1,0){\line(0,1){1}}
\put(2,0){\line(0,1){1}}
\end{picture}\,
}
\lb a_1,a_2,a_3\rb
\\
&=k_3(a_1,a_2,a_3)+
k_1(a_1)k_2(a_2,a_3)
+k_2(a_1,a_2)k_1(a_3)\\
&\quad
+k_2(a_1,a_3)k_1(a_2)
+k_1(a_1)k_1(a_2)k_1(a_3),
\end{align*}
and thus
\begin{align*}
k_3(a_1,a_2,a_3)&=\ff(a_1a_2a_3)-\ff(a_1)\ff(a_2a_3)-\ff(a_1a_3)
\ff(a_2)\\&\quad-\ff(a_1a_2)\ff(a_3)+2\ff(a_1)\ff(a_2)\ff(a_3).
\end{align*}
\end{itemize}
3) For $n=4$ we consider the special case where all
$\ff(a_i)=0$. Then we have
$$k_4(a_1,a_2,a_3,a_4)=\ff(a_1a_2a_3a_4)-\ff(a_1a_2)\ff(a_3a_4)
-\ff(a_1a_4)\ff(a_2a_3).$$
4) The $k_n$ are multi-linear functionals in their
$n$ arguments.

\subsection*{}
The meaning of the concept  `cumulants' for 
freeness is shown by the following theorem.

\subsection{Theorem}
Let $(\cA,\ff)$ be a probability space
and consider unital
subalgebras $\cA_1,\dots,\cA_m\subset\cA$.
Then the following two statements are equivalent:
\begin{itemize}
\item[i)]
$\cA_1,\dots,\cA_m$ are free.
\item[ii)]
We have for all $n\geq 2$ and for all $a_i\in\cA_{j(i)}$ with
$1\leq j(1),\dots,j(n)\leq m$:
$$k_n(a_1,\dots,a_n)=0\qquad
\text{if there exist $1\leq l,k\leq n$ with $j(l)\not=j(k)$.}$$
\end{itemize}

\subsection{Remarks}
1) This characterization of freeness in terms of cumulants is
the translation of the definition of freeness in terms of
moments -- by using the relation between
moments and cumulants from Definition 2.4.
One should note that in contrast to the characterization in
terms of moments we do not require that
$j(1)\not=j(2)\not=\dots\not=j(n)$ or $\ff(a_i)=0$.
Hence the characterization of freeness in terms of
cumulants is much easier to use in concrete calculations.
\\
2) Since the unit $1$ is free from everything, the above
theorem contains as a special case the statement:
$$k_n(a_1,\dots,a_n)=0\qquad\text{if $n\geq 2$ and 
$a_i=1$ for at least one $i$.}$$
This special case will also present an important step in the
proof of Theorem 2.6 and it will be proved 
separately as a lemma. 
\\
3) Note also: for $n=1$  we have
$$k_1(1)=\ff(1)=1.$$

\begin{proof}
(i) $\folgtohne$ (ii): If all $a_i$ 
are centered, i.e. $\ff(a_i)=0$, and 
alternating, i.e.  $j(1)\not=j(2)\not=\dots\not=j(n)$, then the
assertion follows directly by the definition of freeness and
by the relation
$$\ff(a_1\dots a_n)=
k_n(a_1,\dots,a_n)+
\sum_{\pi\in NC(n)\atop \pi\not=\eins_n} k_\pi\lb
a_1,\dots,a_n\rb,$$
because at least one factor of $k_\pi$ for
$\pi\not=\eins_n$ is of the form
$$k_{p+1}(a_l,a_{l+1},\dots,a_{l+p})\qquad\text{with $p+1<n$}$$
and thus the assertion follows by induction.
\\
The essential part of the proof consists in showing that
on the level of cumulants the assumption `centered' is
not needed and `alternating' can be weakened to
`mixed'.
\\
Let us start with getting rid of the assumption `centered'.
For this we will need the following lemma -- which is
of course a special case of our theorem.

\setcounter{subsection}{6}

\subsubsection{Lemma}
Let  $n\geq 2$
und $a_1,\dots,a_n\in\cA$. Then we have:
$$\text{there exists a $1\leq i\leq n$ with $a_i=1$}
\folgt k_n(a_1,\dots,a_n)=0.$$

\begin{proof}
To simplify notation we consider the case
$a_n=1$, i.e. we want to show
$$k_n(a_1,\dots,a_{n-1},1)\zz 0.$$
We will prove this by induction on $n$.
\\
$n=2$ : the assertion is true, since
$$k_2(a,1)=\ff(a1)-\ff(a)\ff(1)=0.$$
$n-1\to n$: Assume we have proved the assertion for all
$k<n$.
Then we have
\begin{align*}
\ff(a_1\dots a_{n-1}1)&=\sum_{\pi\in NC(n)}
k_\pi\lb a_1,\dots,a_{n-1},1\rb\\
&=k_n(a_1,\dots,a_{n-1},1)+\sum_{\pi\in NC(n)\atop
\pi\not=\eins_n} k_\pi\lb a_1,\dots,a_{n-1},1\rb.
\end{align*}
According to our induction hypothesis only such 
$\pi\not=\eins_n$ contribute to the above sum which
have the property that $(n)$ is a one-element block of
$\pi$, i.e. which have the form 
$$\pi=\sigma\cup (n) \qquad\text{with}\qquad
\sigma\in NC(n-1).$$
Then we have
\begin{align*}
k_{\pi}\lb a_1,\dots,a_{n-1},1\rb&=k_{\sigma}\lb
a_1,\dots,a_{n-1}\rb k_1(1)=k_{\sigma}\lb a_1,\dots,a_{n-1}\rb,
\end{align*}
hence
\begin{align*}
\ff(a_1\dots a_{n-1}1)&=k_n(a_1,\dots,a_{n-1},1)+
\sum_{\sigma\in NC(n-1)}k_{\sigma}\lb a_1,\dots,
a_{n-1}\rb\\
&=k_n(a_1,\dots,a_{n-1},1)+\ff(a_1\dots a_{n-1}).
\end{align*}
Since
$$\ff(a_1\dots a_{n-1}1)=\ff(a_1\dots a_{n-1}),$$
we obtain
$$k_n(a_1,\dots,a_{n-1},1)=0. $$
\end{proof}

Let $n\geq 2$. 
Then this lemma implies
that we have for arbitrary
$a_1,\dots,a_n\in\cA$ the relation
$$k_n(a_1,\dots,a_n)=k_n\bigl(a_1-\ff(a_1)1,\dots,a_n
-\ff(a_n)1\bigr),$$
i.e. we can center the arguments of our cumulants $k_n$
($n\geq 2$) without changing the value of the cumulants.
\\
Thus we have proved the following statement:
Consider $n\geq 2$ and $a_i\in\cA_{j(i)}$ ($i=1,\dots,n$) with
$j(1)\not=j(2)\not=\dots\not=j(n)$. 
Then we have
$$k_n(a_1,\dots,a_n)= 0.$$

It remains to weaken the assumption `alternating'  to
`mixed'. For this we will need the following lemma.

\subsubsection{Lemma}
Consider $n\geq 2$, $a_1,\dots,a_n\in\cA$ and
$1\leq p\leq n-1$. Then we have
\begin{multline*}
k_{n-1}(a_1,\dots,a_{p-1},a_pa_{p+1},a_{p+2},\dots,a_n)=
k_n(a_1,\dots,a_p,a_{p+1},\dots,a_n)\\+
\sum_{\pi\in NC(n)\atop \#\pi=2,
p\not\sim_\pi p+1} k_\pi\lb a_1,\dots,a_p,a_{p+1},\dots,
a_n\rb.
\end{multline*}

Examples:
$$k_2(a_1a_2,a_3)=k_3(a_1,a_2,a_3)+k_1(a_1)k_2(a_2,a_3)+
k_2(a_1,a_3)k_1(a_2)$$
\begin{align*}
k_3(a_1,a_2a_3,a_4)&=k_4(a_1,a_2,a_3,a_4)+k_1(a_2)k_3(a_1,a_3,a_4)\\
&\quad
+k_2(a_1,a_2)k_2(a_3,a_4)+k_3(a_1,a_2,a_4)k_1(a_3)
\end{align*}
\begin{proof}
For $\pi\in NC(n)$ we denote by $\pi\vert_{p=p+1}\in
NC(n-1)$ that partition which is obtained by identifying
$p$ and $p+1$, i.e. for $\pi=\{V_1,\dots,V_r\}$
we have 
\begin{multline*}
\pi\vert_{p=p+1}=
\{V_1,\dots,(V_k\cup V_l)\backslash\{p+1\},\dots,V_r\},\qquad
\text{if $p\in V_k$ and $p+1\in V_l$.}
\end{multline*}
(If $p$ and $p+1$ belong to different blocks, then
$\pi\vert_{p=p+1}$ has one block less than $\pi$; 
if $p$ and $p+1$ 
belong to the same block, then the number of blocks does
not change;
of course, we identify partitions of the set $\{1,\dots,p,p+2,\dots,n\}$ 
with partitions from
$NC(n-1)$;  
the property `non-crossing' is preserved under the
transition from $\pi$ to $\pi\vert_{p=p+1}$.)
\\
\setlength{\unitlength}{0.3cm}
Example:
Consider
$$\pi=\{(1,5,6),(2),(3,4)\}\entspricht
\begin{picture}(5,2)\thicklines
\put(0,0){\line(0,1){2}}
\put(0,0){\line(1,0){5}}
\put(4,0){\line(0,1){2}}
\put(5,0){\line(0,1){2}}
\put(1,1){\line(0,1){1}}
\put(2,1){\line(0,1){1}}
\put(2,1){\line(1,0){1}}
\put(3,1){\line(0,1){1}}
\put(-0.3,2.2){1}
\put(0.7,2.2){2}
\put(1.7,2.2){3}
\put(2.7,2.2){4}
\put(3.7,2.2){5}
\put(4.7,2.2){6}
\end{picture}
$$
Then we have
$$\pi\vert_{5=6}=\{(1,5),(2),(3,4)\}\entspricht 
\begin{picture}(4,2)\thicklines
\put(0,0){\line(0,1){2}}
\put(0,0){\line(1,0){4}}
\put(4,0){\line(0,1){2}}
\put(1,1){\line(0,1){1}}
\put(2,1){\line(0,1){1}}
\put(2,1){\line(1,0){1}}
\put(3,1){\line(0,1){1}}
\put(-0.3,2.2){1}
\put(0.7,2.2){2}
\put(1.7,2.2){3}
\put(2.7,2.2){4}
\put(3.7,2.2){5}
\end{picture}
$$
and
$$\pi\vert_{4=5}=\{(1,3,4,6),(2)\}\entsprichtohne
\{(1,3,4,5),(2)\}\entspricht
\begin{picture}(4,2)\thicklines
\put(0,0){\line(0,1){2}}
\put(0,0){\line(1,0){4}}
\put(4,0){\line(0,1){2}}
\put(1,1){\line(0,1){1}}
\put(2,0){\line(0,1){2}}
\put(3,0){\line(0,1){2}}
\put(-0.3,2.2){1}
\put(0.7,2.2){2}
\put(1.7,2.2){3}
\put(2.7,2.2){4}
\put(3.7,2.2){6}
\end{picture}
$$
With the help of this definition we can state our 
assertion more generally for $k_\sigma$
for arbitrary $\sigma\in NC(n-1)$: Assume that our
assertion is true for all $l< n$, i.e.
\begin{align*}
k_{l-1}(a_1,\dots,a_pa_{p+1},\dots,a_l)&=
k_{\eins_{l-1}}\lb a_1,\dots,a_pa_{p+1},\dots,a_l\rb\\&=
\sum_{\pi\in NC(l)\atop \pi\vert_{p=p+1}=\eins_{l-1}}
k_\pi\lb a_1,\dots,a_{p},a_{p+1},\dots,a_l\rb.
\end{align*}
Then it is quite easy to see that 
we have for arbitrary $\sigma\in NC(n-1)$ with $\sigma\not=
\eins_{n-1}$:
$$k_\sigma\lb a_1,\dots,a_pa_{p+1},\dots,a_n\rb
=\sum_{\pi\in NC(n)\atop
\pi\vert_{p=p+1}=\sigma}
k_\pi\lb a_1,\dots,a_p,a_{p+1},\dots,a_n\rb.$$
We will now prove the assertion of our lemma by induction 
$n$.\\
$n=2$: The assertion is true because
\begin{align*}
k_1(a_1a_2)&=\ff(a_1a_2)\\
&=(\ff(a_1a_2)-\ff(a_1)\ff(a_2))+\ff(a_1)\ff(a_2)\\
&=k_2(a_1,a_2)+k_1(a_1)k_1(a_2).
\end{align*}
$n-1\to n$: Let the assertion be proven for all $l<n$,
which implies, as indicated above,
that we have also for all $\sigma\in
NC(n-1)$ with $\sigma\not=\eins_{n-1}$
$$k_\sigma\lb a_1,\dots,a_pa_{p+1},\dots,a_n\rb=
\sum_{\pi\in NC(n)\atop \pi\vert_{p=p+1}=\sigma}
k_\pi\lb a_1,\dots,a_p,a_{p+1},\dots,a_n\rb.$$
Then we have
\begin{align*}
k_{n-1}(a_1&,\dots,a_pa_{p+1},\dots,a_n)\\&=
\ff(a_1\dots(a_pa_{p+1})\dots a_n)-
\sum_{\sigma\in NC(n-1)\atop \sigma\not=\eins_{n-1}}
k_\sigma\lb a_1,\dots,a_pa_{p+1},\dots,a_n\rb \\
&=\ff(a_1\dots a_pa_{p+1}\dots a_n)-
\sum_{\sigma\in NC(n-1)\atop \sigma\not=\eins_{n-1}}
\sum_{\pi\in NC(n)\atop \pi\vert_{p=p+1}=\sigma}
k_\pi\lb a_1,\dots,a_p,a_{p+1},\dots,a_n\rb
\displaybreak[0]
\\
&=\ff(a_1\dots a_pa_{p+1}\dots a_n)-
\sum_{\pi\in NC(n)\atop \pi\vert_{p=p+1}\not=\eins_{n-1}}
k_\pi\lb a_1,\dots,a_p,a_{p+1},\dots,a_n\rb
\displaybreak[0]
\\
&=
\sum_{\pi\in NC(n)}k_\pi\lb a_1,\dots,a_p,a_{p+1},
\dots,a_n\rb-
\sum_{\pi\in NC(n)\atop \pi\vert_{p=p+1}\not=\eins_{n-1}}
k_\pi\lb a_1,\dots,a_p,a_{p+1},\dots,a_n\rb
\displaybreak[0]
\\
&=\sum_{\pi\in NC(n)\atop  \pi\vert_{p=p+1}=\eins_{n-1}}
k_\pi\lb a_1,\dots,a_p,a_{p+1},\dots,a_n\rb\\
&=k_n(a_1,\dots,a_p,a_{p+1},\dots,a_n)+
\sum_{\pi\in NC(n)\atop \#\pi=2, p\not\sim_\pi p+1}
k_\pi\lb a_1,\dots,a_p,a_{p+1},\dots,a_n\rb.
\end{align*}
\end{proof}

By using this lemma we can now prove our theorem in
full generality:
Consider $n\geq 2$ and $a_i\in\cA_{j(i)}$ ($i=1,\dots,n$).
Assume that there exist $k,l$ with $j(k)\not=j(l)$. We have to
show
$$k_n(a_1,\dots,a_n)\zz 0.$$
This follows so:
If $j(1)\not=j(2)\not=\dots\not=j(n)$, 
then the assertion is already proved. Thus we can assume that 
there exists a $p$ with $j(p)=j(p+1)$, implying
$a_pa_{p+1}\in\cA_{j(p)}$. In that case we can use
the above lemma to obtain
\begin{align*}
k_n(a_1,\dots,a_p,a_{p+1},\dots,a_n)&=
k_{n-1}(a_1,\dots,a_pa_{p+1},\dots,a_n)\\&\quad-
\sum_{\pi\in NC(n)\atop \#\pi=2, p\not\sim_\pi p+1}
k_\pi\lb a_1,\dots,a_p,a_{p+1},\dots,a_n\rb.
\end{align*}
To show that this vanishes we will again use induction
on $n$: The first term\linebreak
$k_{n-1}(a_1,\dots,a_pa_{p+1},\dots,a_n)$ 
vanishes by induction hypothesis, since two of its
arguments are lying in the different algebras $\cA_{j(k)}$ and $\cA_{j(l)}$. 
Consider now the summand
$k_\pi\lb a_1,\dots,a_p,a_{p+1},\dots,a_n\rb$
for
$\pi\in NC(n)$ with 
$$\pi=\{V_1,V_2\},\qquad\text{where $p\in V_1$ and 
$p+1\in V_2$.}$$
Then we have
$k_\pi=k_{V_1}k_{V_2}$,
and by induction hypothesis this can be different from zero
only in the case where
all arguments in each of the two factors are coming from
the same algebra; but this would impy that in the first
factor all arguments are in $\cA_{j(p)}$
and in the second factor all arguments are in
$\cA_{j(p+1)}$. Because of
$j(p)=j(p+1)$ this would imply
$j(1)=j(2)=\dots=j(n)$,
yielding a contradiction with
$j(l)\not=j(k)$. 
Thus all terms of the right hand side have to vanish and
we obtain
$$k_n(a_1,\dots,a_p,a_{p+1},\dots,a_n)=0.$$
(ii) $\folgtohne$ (i): 
(ii) gives an inductive way to calculate uniquely all
mixed moments; according to what we have proved
above this mixed moments must calculate in the same way
as for free subalgebras; but this means of course
that these subalgebras are free.
\end{proof}

\setcounter{subsection}{7}

\subsection{Notation}
For a random variable $a\in\cA$ we put
$$k_n^a:=k_n(a,\dots,a)$$
and call $(k_n^a)_{n\geq1}$
the {\bf (free) cumulants of $a$}.

\subsection*{}
Our main theorem on the vanishing of mixed cumulants in
free variables specifies in this one-dimensional case to
the linearity of the cumulants.

\subsection{Proposition}
Let $a$ and $b$ be free. Then we have
$$k_n^{a+b}=k_n^a+k_n^b\qquad\text{for all $n\geq 1$.}$$

\begin{proof}
We have
\begin{align*}
k_n^{a+b}&=k_n(a+b,\dots,a+b)\\
&=k_n(a,\dots,a)+k_n(b,\dots,b)\\
&= k_n^a+k_n^b,
\end{align*}
because cumulants which have both $a$ and $b$ as arguments
vanish by Theorem 2.6.
\end{proof}

\subsection*{}
Thus, free convolution is easy to describe on the level
of cumulants; the cumulants are additive under free
convolution. It remains to make the connection between
moments and cumulants as explicit as possible. On a
combinatorial level, our definition specializes in the 
one-dimensional case to the following relation.

\subsection{Proposition}
Let $(m_n)_{n\geq1}$ 
and $(k_n)_{n\geq 1}$ be the moments and free cumulants,
respectively, of some
random variable. The connection between these two 
sequences of numbers is given by
$$m_n=\sum_{\pi\in NC(n)} k_\pi,$$
where
$$k_\pi:=k_{\#V_1}\cdots k_{\#V_r}\qquad\text{for}\qquad
\pi=\{V_1,\dots,V_r\}.$$

Example: For $n=3$ we have
\setlength{\unitlength}{0.2cm}
\begin{align*}
m_3&=
k_{\,
\begin{picture}(2,1)\thicklines
\put(0,0){\line(0,1){1}}
\put(0,0){\line(1,0){2}}
\put(1,0){\line(0,1){1}}
\put(2,0){\line(0,1){1}}
\end{picture}\,
}
+k_{\,
\begin{picture}(2,1)\thicklines
\put(0,0){\line(0,1){1}}
\put(1,0){\line(1,0){1}}
\put(1,0){\line(0,1){1}}
\put(2,0){\line(0,1){1}}
\end{picture}\,
}
+k_{\,
\begin{picture}(2,1)\thicklines
\put(0,0){\line(0,1){1}}
\put(0,0){\line(1,0){1}}
\put(1,0){\line(0,1){1}}
\put(2,0){\line(0,1){1}}
\end{picture}\,
}
+k_{\,
\begin{picture}(2,1)\thicklines
\put(0,0){\line(0,1){1}}
\put(0,0){\line(1,0){2}}
\put(1,0.5){\line(0,1){0.5}}
\put(2,0){\line(0,1){1}}
\end{picture}\,
}
+k_{\,
\begin{picture}(2,1)\thicklines
\put(0,0){\line(0,1){1}}
\put(1,0){\line(0,1){1}}
\put(2,0){\line(0,1){1}}
\end{picture}\,
}\\
&=k_3+3k_1k_2+k_1^3.
\end{align*}

\subsection*{}
For concrete calculations, however, one would prefer to
have a more analytical description of the relation between
moments and cumulants. This can be achieved by translating
the above relation to corresponding formal power series.

\subsection{Theorem}
Let $(m_n)_{n\geq 1}$ and $(k_n)_{n\geq 1}$ 
be two sequences of complex numbers and consider the
corresponding formal power series 
\begin{align*}
M(z)&:=1+\sum_{n=1}^\infty m_nz^n,\\
C(z)&:=1+\sum_{n=1}^\infty k_nz^n.
\end{align*}
Then the following three statements are equivalent:
\begin{itemize}
\item[(i)]
We have for all $n\in\NN$
$$m_n=\sum_{\pi\in NC(n)}k_\pi=
\sum_{\pi=\{V_1,\dots,V_r\}\in NC(n)}k_{\#V_1}\dots k_{\#V_r}.$$
\item[(ii)]
We have for all $n\in\NN$
(where we put $m_0:=1$)
$$
m_n=\sum_{s=1}^n\sum_{i_1,\dots,i_s\in\{0,1,\dots,n-s\}\atop
i_1+\dots+i_s=n-s} k_sm_{i_1}\dots m_{i_s}.
$$
\item[(iii)]
We have
$$C\lb zM(z)\rb=M(z).$$
\end{itemize}

\begin{proof}
We rewrite the sum
$$m_n=\sum_{\pi\in NC(n)}k_\pi$$
in the way that we fix the first block $V_1$ of $\pi$ (i.e. that
block which contains the element 1) and sum over all possibilities
for the other blocks; in the end we sum over $V_1$:
$$m_n=\sum_{s=1}^n\quad\sum_{\text{$V_1$ with $\#V_1=s$}}
\sum_{\pi\in NC(n)\atop \text{where $\pi=\{V_1,\dots\}$}}
k_\pi.$$
If
$$V_1=(v_1=1,v_2,\dots,v_s),$$
then $\pi=\{V_1,\dots\}\in NC(n)$ 
can only connect elements lying between some
$v_k$ and $v_{k+1}$, i.e.
$\pi=\{V_1,V_2,\dots,V_r\}$ such that we have for all $j=2,\dots,r$:
there exists a $k$ with $v_k<V_j<v_{k+1}$.
There we put
$$v_{s+1}:=n+1.$$
Hence such a $\pi$ decomposes as
$$\pi=V_1\cup \tilde\pi_1\cup\dots\cup\tilde\pi_s,$$
where
$$\text{$\tilde\pi_j$ is a non-crossing partition of
 $\{v_j+1,v_j+2,\dots,
v_{j+1}-1\}$.}$$
For such $\pi$ we have
\begin{align*}
k_\pi&=k_{\#V_1}k_{\tilde \pi_1}\dots k_{\tilde\pi_s}
=k_s k_{\tilde \pi_1}\dots k_{\tilde\pi_s},
\end{align*}
and thus we obtain
\begin{align*}
m_n&=\sum_{s=1}^n\sum_{1=v_1<v_2<\dots<v_s\leq n}
\sum_{\pi=V_1\cup\tilde\pi_1\cup\dots\cup\tilde\pi_s\atop
\tilde\pi_j\in NC(v_j+1,\dots,v_{j+1}-1)}
k_sk_{\tilde\pi_1}\dots k_{\tilde\pi_s}\\
&=\sum_{s=1}^n k_s
\sum_{1=v_1<v_2<\dots<v_s\leq n}
\bigl(\sum_{\tilde\pi_1\in NC(v_1+1,\dots,v_2-1)}
k_{\tilde\pi_1}\bigr)\dots\bigl(
\sum_{\tilde\pi_s\in NC(v_s+1,\dots,n)}
k_{\tilde\pi_s}\bigr)\\
&=\sum_{s=1}^n k_s \sum_{1=v_1<v_2<\dots<v_s\leq n}
m_{v_2-v_1-1}\dots m_{n-v_s}\\
&=\sum_{s=1}^n\sum_{i_1,\dots,i_s\in\{0,1,\dots,n-s\}\atop
i_1+\dots+i_s+s=n} k_sm_{i_1}\dots m_{i_s}
\qquad\qquad\qquad (i_k:=v_{k+1}-v_k-1).
\end{align*}
This yields the implication (i) $\folgtohne$ (ii).
\\
We can now rewrite (ii) 
in terms of the corresponding formal power series in
the following way (where we put $m_0:=k_0:=1$):
\begin{align*}
M(z)&=
1+\sum_{n=1}^\infty z^nm_n \\
&=1+\sum_{n=1}^\infty \sum_{s=1}^n
\sum_{i_1,\dots,i_s\in\{0,1,\dots,n-s\}\atop
i_1+\dots+i_s=n-s} k_sz^s m_{i_1}z^{i_1}\dots m_{i_s}z^{i_s}\\
&=1+\sum_{s=1}^\infty k_s z^s \bigl(\sum_{i=0}^\infty m_iz^i\bigr)^s\\
&=C\lb zM(z)\rb.
\end{align*}
This yields (iii).
\\
Since (iii) describes uniquely a fixed relation between 
the numbers $(k_n)_{n\geq 1}$ and the numbers $(m_n)_{n\geq 1}$, this has to
be the relation (i).
\end{proof}

\subsection*{}
If we rewrite the above relation between the formal
power series in terms of the Cauchy-transform
$$G(z):=\sum_{n=0}^\infty \frac {m_n}{z^{n+1}}$$
and the $R$-transform
$$R(z):=\sum_{n=0}^\infty k_{n+1}z^n,$$
then we obtain Voiculescu's formula.
 
\subsection{Corollary}
The relation between the Cauchy-transform $G(z)$ and the
$R$-transform $R(z)$ of a random variable is given by
$$G\lb R(z)+\frac 1z\rb=z.$$

\begin{proof}
We just have to note that the formal power series $M(z)$ and
$C(z)$ from Theorem 2.11 and $G(z)$, $R(z)$, and 
$K(z)=R(z)+\frac 1z$ 
are related by: 
$$G(z)=\frac 1z M(\frac 1z)$$
and
$$C(z)=1+zR(z)=zK(z),\qquad\text{thus}\qquad
K(z)=\frac {C(z)}z.$$
This gives
\begin{align*}
K\lb G(z)\rb&=\frac 1{G(z)}C\lb G(z)\rb
=\frac 1{G(z)} C\lb \frac 1zM(\frac 1z)\rb=\frac 1{G(z)}M(\frac 1z)=z,
\end{align*}
thus
$K\lb G(z)\rb=z $
and hence also
$$G\lb R(z)+\frac 1z\rb=G\lb K(z)\rb=z.$$
\end{proof}

\subsection{Remark}
It is quite easy to check that the cumulants $k_n^a$ of
a random variable $a$ are indeed the coefficients of
the $R$-transform of $a$ as introduced by 
Voiculescu: Remember that the distribution
of $a$ was modelled by the canonical variable
(special formal power series in an isometry $l^*$,
see \cite{VDN})
$$b=l^*+\sum_{i=0}^\infty k_{i+1} l^i\in(\FPR(l),\tau).$$
Then we have
\begin{align*}
m_n&=\tau\bigl((l^*+\sum_{i=0}^\infty k_{i+1}l^i)^n\bigr)\\
&=\sum_{i(1),\dots,i(n)\in\{-1,0,1,\dots,n-1\}}
\tau(l^{i(n)}\dots l^{i(1)})k_{i(1)+1}\dots k_{i(n)+1},
\end{align*}
where
$l^{-1}$ is identified with $l^*$, 
$$l^{-1}\entspricht l^*$$
and
$$k_0:=1.$$

The sum is running over tuples $(i(1),\dots,i(n))$, 
which can be identified with paths in the lattice
$\ZZ^2$: 
\begin{align*}
i=-1&\entspricht\text{diagonal step upwards:}
\begin{pmatrix}
1\\1
\end{pmatrix}\\
i=0&\entspricht\text{horizontal step to the right:}
\begin{pmatrix}
1\\0
\end{pmatrix}\\
i=k \quad (1\leq k\leq n-1)&\entspricht\text{diagonal step
downwards:}
\begin{pmatrix}
1\\-k
\end{pmatrix}
\end{align*}
We have now
$$\tau(l^{i(n)}\dots  l^{i(1)})=\begin{cases}
1,&\text{if $i(1)+\dots+ i(m)\leq 0$ $\forall\,m=1,\dots,n$ and}\\
&\qquad\text{$i(1)+\dots+ i(n)=0$}\\
0,&\text{otherwise}
\end{cases}$$
and thus
$$m_n=\sum_{{i(1),\dots,i(n)\in\{-1,0,1,\dots,n-1\}\atop
i(1)+\dots +i(m)\leq 0\quad\forall\, m=1,\dots,n}\atop
i(1)+\dots+i(n)=0}
k_{i(1)+1}\dots k_{i(n)+1}.$$

Hence only such paths from $(0,0)$ to $(n,0)$ contribute
which stay always above the $x$-axis.
Each such path is weighted in a multiplicative way 
(using the cumulants) with the
length of its steps.
\\
Example:
\setlength{\unitlength}{0.5cm}
$$
\begin{picture}(6,3)
\thinlines
\multiput(0,0)(0,1){4}{\line(1,0){6}}
\multiput(0,0)(1,0){7}{\line(0,1){3}}
\thicklines
\put(0,0){\vector(1,1){1}}
\put(0,0.5){1}
\put(1,1){\vector(1,1){1}}
\put(1,1.5){1}
\put(2,2){\vector(1,0){1}}
\put(1.8,2.5){$k_1$}
\put(3,2){\vector(1,1){1}}
\put(3,2.5){1}
\put(4,3){\vector(1,-2){1}}
\put(4.5,2.1){$k_3$}
\put(5,1){\vector(1,-1){1}}
\put(6.0,0.5){$k_2$}
\end{picture}
\entspricht
\tau(l^1l^2l^*l^0l^*l^*)k_1k_3k_2
$$

The above summation can now be rewritten in terms of
a summation over non-crossing partitions leading to
the relation from Proposition 2.10. We will leave the proof
of this lemma to the reader 

\subsubsection{Lemma}
There exists a canonical bijection 
\begin{align*}
NC(n)\longleftrightarrow \{(i(1),\dots,i(n))\mid &i(m)\in\{-1,0,1,
\dots,n-1\},\\
&i(1)+\dots+ i(m)\leq 0\quad\forall\, m=1,\dots,n;\\
&i(1)+\dots+ i(n)=0 \qquad\qquad\qquad\qquad\}.
\end{align*}
It is given by
$$\pi\mapsto \Pi=(i(1),\dots,i(n))$$
where
$$i(m)=\begin{cases}
\#V_i-1,&\text{if $m$ is the last element
in a block $V_i$}\\
-1,&\text{otherwise} 
\end{cases}$$

Note that a block consisting of one element gives a horizontal
step; a block consisting of $k$ ($k\geq 2$) elements
gives $k-1$ upward steps each of length 1 and one
downward step of length $k-1$.
\\
An example for this bijection is
$$\setlength{\unitlength}{0.3cm}
\begin{picture}(6,4)\thicklines
\put(0,0){\line(0,1){3}}
\put(0,0){\line(1,0){5}}
\put(5,0){\line(0,1){3}}
\put(1,1){\line(0,1){2}}
\put(1,1){\line(1,0){3}}
\put(3,1){\line(0,1){2}}
\put(4,1){\line(0,1){2}}
\put(2,2){\line(0,1){1}}
\put(-0.3,3.2){1}
\put(0.7,3.2){2}
\put(1.7,3.2){3}
\put(2.7,3.2){4}
\put(3.7,3.2){5}
\put(4.7,3.2){6}
\end{picture}
\qquad\entsprichtohne\quad
\pi=\{(1,6),(2,4,5),(3)\}$$
is mapped to
$$\Pi=(-1,-1,0,-1,2,1)
\quad\entsprichtohne\qquad
\setlength{\unitlength}{0.4cm}
\begin{picture}(6,3)
\thinlines
\multiput(0,0)(1,0){7}{\line(0,1){3}}
\multiput(0,0)(0,1){4}{\line(1,0){6}}
\thicklines
\put(0,0){\vector(1,1){1}}
\put(1,1){\vector(1,1){1}}
\put(2,2){\vector(1,0){1}}
\put(3,2){\vector(1,1){1}}
\put(4,3){\vector(1,-2){1}}
\put(5,1){\vector(1,-1){1}}
\end{picture}
$$

Now note that with this identification of paths and
non-crossing partitions the factor
$$k_\Pi=k_{i(1)+1}\dots k_{i(n)+1}$$
for
$$\Pi=(i(1),\dots,i(n))
\entspricht\pi=\{V_1,\dots,V_r\}$$
goes over to
$$k_\pi:=k_{\#V_1}\dots k_{\#V_r}.$$
Consider the above example:
$$\pi=\{(1,6),(2,4,5),(3)\}\mapsto\Pi\quad
\entsprichtohne\qquad
\setlength{\unitlength}{0.4cm}
\begin{picture}(6,3)
\thinlines
\multiput(0,0)(1,0){7}{\line(0,1){3}}
\multiput(0,0)(0,1){4}{\line(1,0){6}}
\thicklines
\put(0,0){\vector(1,1){1}}
\put(1,1){\vector(1,1){1}}
\put(2,2){\vector(1,0){1}}
\put(3,2){\vector(1,1){1}}
\put(4,3){\vector(1,-2){1}}
\put(5,1){\vector(1,-1){1}}
\put(1.8,2.5){$k_1$}
\put(4.5,2.1){$k_3$}
\put(6.0,0.5){$k_2$}
\end{picture}
$$
thus
$$k_\pi=k_\Pi=k_1k_3k_2=k_{\#(3)}
k_{\#(2,4,5)}k_{\#(1,6)}.$$

This correspondence leads of course to the relation
as stated in Proposition 2.10 ; thus the coefficients of the 
$R$-transform of Voiculescu coincide indeed with the
free cumulants as defined in 2.8.
Note that in this way we
obtained easy combinatorial proofs of the main
facts on the $R$-transform -- namely, its additivity
under free convolution and the formula relating it
to the Cauchy-transform.

\subsection*{}
Finally, to show that our description of freeness in terms
of cumulants has also a significance apart from dealing
with additive free convolution, we will apply it to the
problem of the product of free random variables:
Consider
$a_1,\dots,a_n,b_1,\dots,b_n$ such that
$\{a_1,\dots,a_n\}$ and $\{b_1,\dots,b_n\}$ are free.
We want to express the distribution of the random
variables 
$a_1b_1,\dots,a_nb_n$
in terms of the distribution of the $a$'s and of 
the $b$'s.

\subsection{Notation}
1) Analogously to $k_\pi$ 
we define for 
$$\pi=\{V_1,\dots,V_r\}\in NC(n)$$
the expression
$$\ff_\pi\lb a_1\dots,a_n\rb:=\ff_{V_1}\lb a_1,\dots,a_n\rb
\dots \ff_{V_r}\lb a_1,\dots,a_n\rb,$$
where
$$\ff_V\lb a_1,\dots,a_n\rb:=\ff(a_{v_1}\dots a_{v_l})\qquad
\text{for}\qquad V=(v_1,\dots,v_l).$$
Examples:
\setlength{\unitlength}{0.2cm}
\begin{align*}
\ff_{\,
\begin{picture}(2,1)\thicklines
\put(0,0){\line(0,1){1}}
\put(0,0){\line(1,0){2}}
\put(1,0){\line(0,1){1}}
\put(2,0){\line(0,1){1}}
\end{picture}\,
}
\lb a_1,a_2,a_3\rb
&=\ff(a_1a_2a_3)\\
\ff_{\,
\begin{picture}(2,1)\thicklines
\put(0,0){\line(0,1){1}}
\put(1,0){\line(1,0){1}}
\put(1,0){\line(0,1){1}}
\put(2,0){\line(0,1){1}}
\end{picture}\,
}
\lb a_1,a_2,a_3\rb
&=\ff(a_1)\ff(a_2a_3)\\
\ff_{\,
\begin{picture}(2,1)\thicklines
\put(0,0){\line(0,1){1}}
\put(0,0){\line(1,0){1}}
\put(1,0){\line(0,1){1}}
\put(2,0){\line(0,1){1}}
\end{picture}\,
}
\lb a_1,a_2,a_3\rb
&=\ff(a_1a_2)\ff(a_3)\\
\ff_{\,
\begin{picture}(2,1)\thicklines
\put(0,0){\line(0,1){1}}
\put(0,0){\line(1,0){2}}
\put(1,0.5){\line(0,1){0.5}}
\put(2,0){\line(0,1){1}}
\end{picture}\,
}
\lb a_1,a_2,a_3\rb
&=\ff(a_1a_3)\ff(a_2)\\
\ff_{\,
\begin{picture}(2,1)\thicklines
\put(0,0){\line(0,1){1}}
\put(1,0){\line(0,1){1}}
\put(2,0){\line(0,1){1}}
\end{picture}\,
}
\lb a_1,a_2,a_3\rb
&=\ff(a_1)\ff(a_2)\ff(a_3)
\end{align*}
2) Let $\sigma,\pi\in NC(n)$.
Then we write
$$\sigma\leq \pi$$
if each block of
$\sigma$ is contained as a whole in some 
block of $\pi$, i.e. $\sigma$ can be obtained out of  $\pi$ by
refinement of the block structure.
\\
Example:
$$\{(1),(2,4),(3),(5,6)\}\leq \{(1,5,6),(2,3,4)\}$$

\subsection*{}
With these notations we can generalize the relation
$$\ff(a_1\dots a_n)=\sum_{\pi\in NC(n)}
k_\pi\lb a_1,\dots,a_n\rb$$
in the following way.

\subsection{Proposition}
Consider $n\in\NN$, $\sigma\in NC(n)$ 
and $a_1,\dots,a_n\in\cA$. Then we have 
$$\ff_\sigma\lb a_1,\dots,a_n\rb=
\sum_{\pi\in NC(n)\atop \pi\leq \sigma} k_\pi
\lb a_1,\dots,a_n\rb.$$

\begin{proof}
Each
$$\pi\leq\sigma=\{V_1,\dots,V_r\}$$
can be decomposed as
$$\pi=\pi_1\cup\dots\cup \pi_r\qquad\text{where}\qquad
\pi_i\in NC(V_i)\quad (i=1,\dots,r).$$
In such a case we have of course
$$k_\pi=k_{\pi_1}\dots k_{\pi_r}.$$
Thus we obtain (omitting the arguments)
\begin{align*}
\ff_\sigma&=\ff_{V_1}\dots\ff_{V_r}\\
&=\big(\sum_{\pi_1\in NC(V_1)}k_{\pi_1}\bigr)\dots
\big(\sum_{\pi_r\in NC(V_r)}k_{\pi_r}\bigr)\\
&=\sum_{\pi=\pi_1\cup\dots\cup\pi_r\leq\sigma}k_{\pi_1}
\dots k_{\pi_r}\\
&=\sum_{\pi\leq\sigma}k_\pi.
\end{align*}
\end{proof}

\subsection*{}
Consider now
$$\{a_1,\dots,a_n\},\{b_1,\dots,b_n\}\qquad\text{free}.$$
We want to express alternating moments in $a$ and $b$
in terms of moments of $a$ and moments of $b$.
We have
\begin{align*}
\ff(a_1b_1a_2b_2\dots a_nb_n)&=\sum_{\pi\in NC(2n)}
k_\pi\lb a_1,b_1,a_2,b_2,\dots,a_n,b_n\rb.
\end{align*}
Since the $a$'s 
are free from the $b$'s, Theorem 2.6 tells us that only such
$\pi$ contribute to the sum whose blocks do not connect
$a$'s with $b$'s.
But this means that such a $\pi$ has to decompose as 
\begin{align*}
\pi=\pi_1\cup\pi_2\qquad\text{where}\quad&
\pi_1\in NC(1,3,5,\dots,2n-1)\\
&\pi_2\in NC(2,4,6,\dots,2n).
\end{align*}
Thus we have
\begin{align*}
\ff(a_1b_1a_2b_2\dots a_nb_n)&=\sum_{\pi_1\in 
NC(odd),
\pi_2\in NC(even)\atop \pi_1\cup\pi_2\in NC(2n)}
k_{\pi_1}\lb a_1,a_2,\dots,a_n\rb\cdot
k_{\pi_2}\lb b_1,b_2,\dots,b_n\rb\\
&=\sum_{\pi_1\in NC(odd)}
\Bigl(k_{\pi_1}\lb a_1,a_2,\dots,a_n\rb\cdot
\sum_{\pi_2\in NC(even)\atop
\pi_1\cup\pi_2\in NC(2n)} 
k_{\pi_2}\lb b_1,b_2,\dots,b_n\rb\Bigr).
\end{align*}

Note now that for a fixed $\pi_1$ there exists a maximal element
$\sigma$ with the property
$\pi_1\cup\sigma\in NC(2n)$ and that the second sum
is running over all $\pi_2\leq \sigma$. 

\subsection{Definition}
Let $\pi\in NC(n)$  
be a non-crossing partition of the numbers
$1,\dots,n$. Introduce additional numbers
$\bar 1,\dots,\bar n$, with alternating order between the
old and the new ones, i.e. we order them in the way
$$1\bar 1 2 \bar 2\dots n\bar n.$$
We define the 
{\bf complement $K(\pi)$} of $\pi$
as the maximal  
$\sigma\in NC(\bar 1,\dots,\bar n)$ with the property
$$\pi\cup\sigma\in NC(1,\bar1,\dots,n,\bar n).$$
If we present the partition 
$\pi$ graphically by connecting the blocks in
$1,\dots,n$, then $\sigma$ is given by 
connecting as much as possible the numbers $\bar 1,\dots, \bar n$
without getting crossings among themselves and with $\pi$.

(This natural notation of the complement of a
non-crossing partition is also due to Kreweras \cite{Kr}.
Note that there is no analogue of this for the case of all partitions.)

\subsection*{}
With this definition we can continue our above calculation as
follows:
\begin{align*}
\ff(a_1b_1a_2b_2\dots a_nb_n)
&=\sum_{\pi_1\in NC(n)}
\Bigl(k_{\pi_1}\lb a_1,a_2,\dots,a_n\rb\cdot
\sum_{\pi_2\in NC(n)\atop
\pi_2\leq K(\pi_1)} 
k_{\pi_2}\lb b_1,b_2,\dots,b_n\rb\Bigr)\\
&=\sum_{\pi_1\in NC(n)}
k_{\pi_1}\lb a_1,a_2,\dots,a_n\rb\cdot
\ff_{K(\pi_1)}\lb b_1,b_2,\dots,b_n\rb.
\end{align*}
Thus we have proved the following result.

\subsection{Theorem}
Consider
$$\{a_1,\dots,a_n\},\{b_1,\dots,b_n\}\qquad\text{free}.$$
Then we have
$$\ff(a_1b_1a_2b_2\dots a_nb_n)=
\sum_{\pi\in NC(n)}
k_{\pi}\lb a_1,a_2,\dots,a_n\rb\cdot
\ff_{K(\pi)}\lb b_1,b_2,\dots,b_n\rb.$$

\subsection*{}
Examples: For $n=1$ we get
\begin{align*}
\ff(ab)&=k_1(a)\ff(b)=\ff(a)\ff(b);
\end{align*}
$n=2$ yields
\begin{align*}
\ff(a_1b_1a_2b_2)&=k_1(a_1)k_1(a_2)\ff(b_1b_2)+
k_2(a_1,a_2)\ff(b_1)\ff(b_2)\\
&=\ff(a_1)\ff(a_2)\ff(b_1b_2)+\bigl(\ff(a_1a_2)-\ff(a_1)\ff(a_2)\bigr)
\ff(b_1)\ff(b_2)\\
&=\ff(a_1)\ff(a_2)\ff(b_1b_2)+\ff(a_1a_2)\ff(b_1)\ff(b_2)
-\ff(a_1)\ff(a_2)\ff(b_1)\ff(b_2).
\end{align*}

\section{{\bf Free stochastic calculus}}

\subsection*{}
In this lecture, we will develop the analogue of a stochastic
calculus for free Brownian motion. Free Brownian motion
is characterized by the same requirements as classical 
Brownian motion, one only has to replace `independent
increments' by `free increments' and the normal distribution
by the semi-circle. In the same way as classical Brownian 
motion can be written as $a_t+a^*_t$ for $a_t$ and $a^*_t$ being
annihilation and creation operators, respectively, on the
Bosonic Fock space, the free Brownian motion has a canonical
realization as $l_t+l^*_t$ for $l_t$ and $l^*_t$ being (left) annihilation
and creation operators on the full Fock space.
Thus, instead of developping a stochastic calculus for
free Brownian motion $S_t=l_t+l^*_t$, one could also
split $S_t$ into its two summands and develop
a free stochastic calculus for $l_t$ and $l^*_t$, in analogy
to the Hudson-Parthasarathy calculus for $a_t$ and
$a^*_t$. This was done by K\"ummerer and Speicher 
\cite{KSp}. The free stochastic calculus with
respect to $S_t$, which is due to Biane and Speicher
\cite{BSp}, however, has some advantages and we will
here restrict to that theory.

In our presentation we will put the emphasis on two main points:
\begin{itemize}
\item
{\bf appropriate norms:} on a linear level all stochastic calculi
have formally the same structure, the main point lies in establishing
the integrals with respect to appropriate norms; in contrast to
all other known examples, the free calculus has the very strong
feature that one has estimates with respect to the uniform operator
norm; i.e. the free stochastic integrals can be defined in
$L^p$ with $p=\infty$
\item
{\bf Ito formula:} on a formal level the difference between 
stochastic calculi lies in their multiplicative structure; at
least formally, a stochastic calculus is characterized by its
Ito formula; for free stochastic calculus this is very similar
to the Ito-formula for classical Brownian motion, however, due
to non-commutativity there is a small, but decisive difference
\end{itemize}

\subsection{Definition}
A {\bf free Brownian motion} consists of 
\begin{itemize}
\item
a von Neumann algebra
$\cA$ 
\item
a faithful normal tracial state $\tau$ on $\cA$
\item
a filtration $(\cA_t)_{t\geq 0}$  -- i.e. $\cA_t$ are 
von Neumann subalgebras of $\cA$ with 
$$\cA_s\subset \cA_t\qquad\text{for $s\leq t$.}$$
\item
a family of operators $(S_t)_{t\geq 0}$ with
\begin{itemize}
\item
$S_t=S_t^*\in\cA_t$
\item
for each $t\geq 0$, the distribution of $S_t$
is a semi-circle with mean 0 and variance $t$
\item
for all $0\leq s<t$,
the increment $S_t-S_s$ is free from $\cA_s$ 
\item
for all $0\leq s<t$, the distribution of 
the increment $S_t-S_s$ is a semi-circle with mean 0 and
variance $t-s$
\end{itemize}
\end{itemize}
Usually, we will call $(S_t)_{t\geq 0}$ the free Brownian motion.

\subsection*{}
In the same way as the classical Brownian motion can be
realized on Bosonic Fock space, the free Brownian motion
has a concrete realization on the full Fock space
-- as follows by Proposition 1.3. Note, however,
that for the development of our free stochastic calculus we will
not need this concrete realization but just the abstract
properties of $(S_t)_{t\geq 0}$.

\subsection{Theorem}
Let 
$$l_t:=l(\chi_{(0.t)}),\qquad l_t^*=l^*(\chi_{(0.t)})$$
be the left annihilation and creation operators
for the characteristic functions of the interval $(0,t)$
on the full Fock space $\cF(\HH)$ for $\HH=L^2(0,\infty)$.
Put
$$\tau\lb A \rb:=\la \Omega,A\Omega\ra,$$
and
$$S_t:=l_t+l_t^*.$$
Then $(S_t)_{t\geq 0}$ is a free Brownian motion with
respect to the filtration
$$\cA_t:=vN(S_s\mid s\leq t).$$

\subsection{Remark}
According to the connection between freeness and random matrices
there is also a random matrix realization of free Brownian
motion:

Consider random matrices
$$B_t^{(N)}:=\frac 1{\sqrt N} \bigl(B_{ij}(t)\bigr)_{i,j=1}^N,$$
where
\begin{itemize}
\item
$B_{ij}(t)$ are classical real-valued Brownian motions for all $i,j$
\item
the matrices $B_t^{(N)}$ are symmetric, i.e. $B_{ij}(t)=B_{ji}(t)$
for all $i,j$
\item
apart from the symmetry condition, all entries are independent, i.e.
$\{B_{ij}(\cdot)\mid 1\leq i\leq j<\infty\}$ are independent
Brownian motions.
\end{itemize}
Consider now the state 
$$\ff:=\text{E}\circ (\frac 1N \text{tr}),$$
where $\text{E}$ denotes the expectation with respect to the
above specified probability space and $\frac 1N\text{tr}$ is
the normalized trace on $N\times N$ matrices.
\\
Then we have
$$S_t\entspricht \lim_{N\to\infty} B_t^{(N)},$$
i.e. for all $n\in\NN$ and $t_1,\dots,t_n\geq 0$ we have:
$$\tau\lb S_{t_1}\dots S_{t_n}\rb=
\lim_{N\to\infty}\text{E}\lb\frac 1N\text{tr}(B_{t_1}^{(N)}\dots
B_{t_n}^{(N)})\rb.$$

We will not use this realization, but it shows that
our free stochastic calculus can also be viewed as the 
large $N$ limit of stochastic calculus with respect to 
$N\times N$ hermitian matrix valued Brownian motion.

\subsection{Remarks} 
1) Note that we have the non-commutative $L^p$-spaces
associated with our free Brownian motion. Namely,
$L^p(\cA$), for $1\leq p<\infty$, is the completion of $\cA$
with respect to the norm
$$\Vert A\Vert_{L^p}:=\tau\lb\vert A\vert^p\rb^{1/p}.$$
For $p=\infty$, we put 
$$\Vert A\Vert_{L^\infty}:=\Vert A\Vert,\qquad\text{i.e.}\qquad
L^\infty(\cA)=\cA.$$
In the concrete realisation of $S_t$ on the full Fock space,
we can identify $L^2(\cA)$ with the full Fock space $\cF(\HH)$
and we can embed $\cA$ into the full Fock space by the
injective mapping
\begin{align*}
\cA&\subset \cF(\HH)\\
A&\mapsto A\Omega.
\end{align*}
2) For our latter norm estimates it will be important that we
can obtain the operator norm as the limit $p\to\infty$ of
the $L^p$-norms:
For $A\in\cA$ one has
\begin{align*}
\Vert A\Vert&=\lim_{p\to\infty}\Vert A\Vert_{L^p}
=\lim_{m\to\infty}\tau\bigl\lb (A^*A)^m\bigr\rb^{1/2m}.
\end{align*}

\subsection{Remarks}
1) Let $A_t,B_t$ be adapted processes. Then we want to define
the stochastic integral
$$\int A_tdS_tB_t.$$
In contrast to the stochastic theories considered in the other
courses, we
have to face now the new phenomenon of two-sided integrals.
In the usual cases, adaptedness of the process implies that
the differentials commute with the process (or anti-commute in
the fermionic case), thus a two-sided integral can always be reduced
to a one-sided one and there is no need to consider two-sided
integrals. But in our case there is no such reduction. Adaptedness
implies that the differential and the process are free, but this
does not result in any commutation relation.
Thus we should consider as the most general integral the
two-sided one. Note that one could of course restrict to one-sided
integrals in the beginning, 
but then a meaningful form of Ito 
formulas would result automatically in two-sided integrals.
\\
2) Since $\int A_tdS_tB_t$ is bilinear in $A_t,B_t$ it is natural
to consider more general linear combinations
$$\sum_i\int A_t^idS_tB_t^i$$
for adapted processes $A_t^i$ and $B_t^i$.
We will also write this as
$$\int U_t\sharp dS_t,\qquad\text{with}\qquad
U_t:=\sum_i A_t^i\otimes B_t^i\in\cA\otimes \cA^{op}$$
and call $U=(U_t)_{t\geq 0}$ a {\bf biprocess}.
($\cA^{op}$ is the opposite algebra of $\cA$, i.e. with the same
linear structure and the order of multiplication reversed; it is quite
natural to consider $U$ as an element in this space, since $A_t$
multiplies from the left, wheras $B_t$ multiplies from the right.)
\\
3) The definition of the integral proceeds now as usual:
First define the integral for simple biprocesses, prove some
adequate norm estimates for such cases and then extend the
definition with respect to the involved norms.

\subsection{Definition}
Let $U_t=A_t\otimes B_t$ be a simple adapted biprocess, i.e.
there exist $0=t_0<t_1<\dots < t_n<\infty$ such that
$$U_t=
\begin{cases}
A_i\otimes B_i& t_i\leq t<t_{i+1}\\
0& t_n\leq t.
\end{cases}$$
(Adaptedness means here of course: $A_i,B_i\in\cA_{t_i}$.)
For such a simple biprocess we define the integral
$$\int U_t\sharp dS_t=\int A_tdS_tB_t:=
\sum_{i=0}^{n-1}A_i(S_{t_{i+1}}-S_{t_i})B_i.$$
For simple adapted biprocesses of the
general form $U_t=\sum_i A_t^i\otimes B_t^i$
we extend the definition by linearity.

\subsection*{}
As usual, it is quite simple to obtain the isometry
of the integral in $L^2$-norm.

\subsection{Proposition (Ito isometry)}
For all adapted simple biprocesses $U$ and $V$, one
has
$$\tau\bigl\lb\int U_t\sharp dS_t\cdot(\int V_t\sharp dS_t)^*
\bigr\rb=\la U,V\ra:=\int \la U_t,V_t\ra_{L^2\otimes L^2}dt.$$

\begin{proof}
By bilinearity, it is enough to prove the assertion for
processes
$$U_t=A\otimes B \cdot 1_{\lb t_0,t_1\lb}(t)\qquad
\text{and}\qquad
V_t=C\otimes D \cdot 1_{\lb t_2,t_3\lb}(t).$$
Then the left hand side is
$$\tau\bigl\lb A(S_{t_1}-S_{t_0})BD^*(S_{t_3}-S_{t_2})C^*\bigr\rb.$$
Note that by linearity it suffices to consider the cases
where the two time intervals are either the same or
disjoint. In the first case we have
\begin{align*}
\tau\bigl\lb A(S_{t_1}-S_{t_0})BD^*(S_{t_1}-S_{t_0})C^*\bigr\rb&=
\tau\lb AC^*\rb\tau\lb BD^*\rb\tau\bigl\lb
(S_{t_1}-S_{t_0})^2\bigr\rb\\
&=\tau\lb AC^*\rb\tau\lb BD^*\rb(t_1-t_0)
\end{align*}
(because $\{A,B,D^*,C^*\}$ is free from the increment
$S_{t_1}-S_{t_0}$), whereas in the second case one of the
increments is free from the rest and thus, because of
the vanishing mean of the increment, the expression
vanishes.
But his gives exactly the assertion.
\end{proof}

\subsection{Notation}
We endow the vector space of simple biprocesses with the
norms ($1\leq p\leq \infty$)
$$\Vert U\Vert_{\cB_p}:=\Bigl(\int\Vert U_t\Vert^2_{
L^p(\tau\otimes\tau^{op})}dt\Bigr)^{1/2}.$$
(For $p=\infty$, $L^\infty(\tau\otimes\tau^{op})$ is the 
von Neumann algebra tensor product of $\cA$ and
$\cA^{op}$.)
\\
The completion of the space of simple biprocesses for these
norms will be denoted by $\cB_p$. The closed subspaces
of $\cB_p$ generated by adapted simple processes will
be denoted by $\cB_p^a$.

\subsection*{}
Thus we have shown that the map $U\mapsto \int U_t\sharp
dS_t$ can be extended isometrically to a mapping
$$\cB_2^a\to L^2(\cA).$$
But it is now an essential feature of the free calculus,
which distinguishes it from all other ones, that one can even
show a norm estimate for $p=\infty$, i.e. the stochastic
integral is for a quite big class of (bi)processes a bounded
operator. Thus we do not have to think about possible
domains of our operators and the multiplication of such
integrals (as considered for the Ito formula) presents no
problems.

\subsection{Theorem (Burkholder-Gundy inequality)}
For any simple adapted biprocess $U$ one has
$$\Vert \int U_t\sharp dS_t\Vert\leq 2\sqrt 2\Vert
U\Vert_{\cB_\infty}.$$

\subsection{Corollary}
The stochastic integral map
$U\mapsto \int U_t\sharp dS_t$ can be extended 
continuously to a mapping
$$\cB_\infty^a\to \cA.$$
In particular, the stochastic integral of an adapted bounded
biprocess from $\cB_\infty^a$ is a bounded operator.

\begin{proof}
Let us just give a sketch of the proof. We restrict here
to biprocesses of the form $U_t=A_t\otimes B_t$. The extension
to sums of such biprocesses follows the same ideas.
\\
Put
$$M_s:=\int_0^sU_t\sharp dS_t=\int_0^s A_tdS_tB_t.$$
We want to obtain an operator norm estimate for $M_s$
by using
$$\Vert M_s\Vert =\lim_{m\to\infty}\tau\bigl\lb(M_s^*M_s)^m
\bigr\rb^{1/2m}.$$
This means we must estimate the $p$-th moment of our
integral for $p\to\infty$. This is much harder than the
case $p=2$, but nevertheless it can be done by using again
the crucial property
$$\tau\bigl\lb A(S_{t_1}-S_{t_0})B(S_{t_1}-S_{t_0})C\bigr\rb=
\tau\lb AC\rb\tau\lb B\rb(t_1-t_0),$$
if $\{A,B,C\}$ is free from the increment $S_{t_1}-S_{t_0}$.
By using also H\"older inequality for non-commutative
$L^p$-spaces one can finally derive an inequality of the form
$$\tau\lb\vert M_s\vert^{2m}\rb\leq 2m\sum_{k=0}^{m-1}
\int_0^s\bigl(\tau\lb\vert M_t\vert^{2k}\rb\cdot
\tau\lb\vert M_t\vert^{2m-2-2k}\rb\cdot \Vert A_t\Vert^2
\cdot \Vert B_t\Vert^2\bigr)dt.$$
Note that the structure of this inequality resembles the
recursion formula for the Catalan numbers $c_n$,
$$c_n=\sum_{k=0}^{n-1}c_kc_{n-k-1}.$$
By induction, one derives now from the above implicit inequality
the explicit one
$$\tau\lb\vert M_s\vert^{2m}\rb\leq c_m\bigl(2\int_0^s\Vert
A_t\Vert^2 \Vert B_t\Vert^2 dt\bigr)^m.$$
We take now the $2m$-th root and note that
$$\lim_{m\to\infty}c_m^{1/2m}=2.$$
Thus we obtain
\begin{align*}
\Vert M_s\Vert &\leq 2 \sqrt 2\bigl(\int_0^s\Vert A_t\Vert^2
\Vert B_t\Vert^2 dt\bigr)^{1/2}\\
&=2\sqrt2 \bigl(\int_0^s \Vert U_t\Vert^2_{L^\infty(\tau\otimes
\tau^{op})}dt\bigr)^{1/2}.
\end{align*}
$s=\infty$ gives the assertion.
\end{proof}

\subsection*{}
Let us now present the Ito formula for the free calculus.
First, we will do this on a formal differential level.
As stated above, if we work with biprocesses from $\cB_\infty^a$,
then our integrals are bounded operators and multiplication 
presents no problem. We will show that the Ito formula holds
even with respect to operator norm.
\\
On a formal level, the Ito formula makes the difference between
different stochastic calculi. 
On a first look, free Brownian motion $(S_t)_{t\geq 0}$ 
has the same Ito formula
as classical Brownian motion $(B_t)_{t\geq 0}$, namely
$$dB_tdB_t=dt \qquad\text{and}\qquad dS_tdS_t=dt.$$
In our non-commutative context, however, this does not contain
all necessary information, since we must now also specify
$$dS_tAdS_t\qquad\text{for $A\in\cA_t$}.$$
In the classical case, $A$ commutes with the increment $dB_t$ and
we have there
$$dB_tAdB_t=AdB_tdB_t=Adt.$$
But for free Brownian motion we have a different result.

\subsection{Theorem (Ito formula -- product form)}
For a free Brownian motion $(S_t)_{t\geq 0}$ we have
the following Ito formula
$$dS_tAdS_t=\tau\lb A\rb dt\qquad\text{for $A\in\cA_t$}.$$

\begin{proof}
Let $I\subset \lb 0,\infty)$ be an interval and consider
decompositions into disjoint sub\-intervals $I_i$,
$I=\cup I_i$. For an interval $I$ we denote 
by $S(I)$ the corresponding increment of the free Brownian motion, i.e.
$$S(I):=S_t-S_s \qquad\text{for}\qquad I=\lb s,t\lb.$$
The main point is now to show that (with $\lambda$ denoting 
Lebesgue measure) 
$$\sum_{i} S(I_i)AS(I_i)
\to\tau\lb A\rb \lambda(I),$$
where we take the usual limit with width $\max_i\lambda(I_i)$
of our decomposition going to zero.
As said above we want to see that this convergence even holds
in operator norm.
\\
We will sketch two proofs of this
fact, one using the abstract properties of
freeness, whereas the other works in a concrete representation
on full Fock space.
\\
1) The assertion follows from the following two facts about freeness:
\begin{itemize}
\item
let $\{s_1,\dots,s_n\}$ be a semicircular family, i.e. each $s_i$
is semicircular and $s_1,\dots,s_n$ are free; then, for a random
variable $a$ which is free from $\{s_1,\dots,s_n\}$ we have
\cite{NS2}:
$s_1as_1,\dots,s_nas_n$ are free
\item
let $x_1,\dots,x_n$ be free random variables with $\tau\lb x_i\rb=0$
for all $i=1,\dots,n$; then we have \cite{V2}
$$\Vert x_1+\dots+x_n\Vert\leq \max_{1\leq i\leq n} \Vert x_i\Vert
+2\bigl(\sum_{i=1}^n\tau\lb\vert x_i\vert^2\rb\bigr)^{1/2}$$
\end{itemize}
2) We realize $S(I)$ on the full Fock space as $S(I)=l(I)+l^*(I)$;
then we have to estimate the four terms
$\sum l(I_i)Al(I_i)$, $\sum l(I_i)Al^*(I_i)$, $\sum l^*(I_i)Al(I_i)$,
and $\sum l^*(I_i)Al^*(I_i)$.
Three of these terms tend to zero by simple norm estimates, the
only problematic case is 
$$\sum_i l(I_i)Al^*(I_i)\to\tau\lb A\rb \lambda(I).$$
(This corresponds of course to the Ito formulas for $l_t$ and $l_t^*$,
namely the only non-zero term is
$dl_tAdl_t^*=\tau\lb A\rb dt$, see \cite{KSp}.)
To prove this later statement, one can model $A$ by the sum of
creation and annihilation operators on the full Fock space as
$A=\sum\alpha_n (l+l^*)^n$ for $l=l(f)$ with $f$ being orthogonal to
$L^2(I)$, hence $l,l^*$ free from all $l(I_i),l^*(I_i)$.
Then one has to expand this representation of $A$ and bring
it, by using the Cuntz relations $l(f)l^*(g)=\la f,g\ra 1$, 
into a normal ordered form
$A=\sum \beta_{n,m} l^{*n}l^m$.
Finally, note that, again by the Cuntz relations, 
only the term for $n=m=0$ contributes
to the sum of our statement.
\end{proof}

\subsection{Example}
The Ito formula contains also the germ for the combinatorial
difference between independence and freeness --
all partitions versus non-crossing partitions.
For a classical Brownian motion $(B_t)_{t\geq 0}$ we can calculate the
fourth moment $\tau\lb B_t^4\rb$ with the help of the Ito formula
as follows:
$$d(B_t^4)=3B_t^3dB_t+6B_t^2dt,$$
thus
$$\frac {d\tau\lb B_t^4\rb}{dt}=6\tau\lb B_t^2\rb=6t$$
yielding
$$\tau\lb B_t^4\rb=3t^2.$$
In the case of the free Brownian motion we obtain in the
same way
\begin{align*}
d(S_t^4)&=dS_tS_t^3+S_tdS_tS_t^2+S_t^2dS_tS_t+S_t^3dS_t\\
&\quad+dS_tdS_tS_tS_t+dS_tS_tdS_tS_t+dS_tS_tS_tdS_t\\
&\quad+S_tdS_tdS_tS_t+
S_tdS_tS_tdS_t+S_tS_tdS_tdS_t\\
&=dS_tS_t^3+S_tdS_tS_t^2+S_t^2dS_tS_t+S_t^3dS_t
+ 3S_t^2dt+\tau\lb S_t^2\rb dt,
\end{align*}
thus
$$\frac{d\tau\lb S_t^4\rb}{dt}=4\tau\lb S_t^2\rb =4t$$
yielding
$$\tau\lb S_t^4\rb=2t^2.$$
This difference between the fourth moments in the classical
and free case reflects the fact that there are 3 pairings of
4 elements, but only 2 of them are non-crossing.

\subsection*{}
The Ito formula can also be put into a functional form to
calculate the differential of a function $f(S_t)$ for nice
functions $f$ -- not just for polynomials.
The main message of the classical Ito formula is that we have to
make a Taylor expansion, but we should not stop after
the first order in the differentials, but take also the second
order into account using $dB_tdB_t=dt$, hence
\begin{align*}
df(B_t)&=f'(B_t)dB_t+\frac 12 f''(B_t)dB_tdB_t\\
&=f'(B_t)dB_t+\frac 12 f''(B_t)dt
\end{align*}
There exists also a free analogue of this; whereas the first order
term is essentially the same as in the classical case, the second
order term is different; one of the two derivatives is replaced
by a difference expression.

\subsection{Theorem (Ito formula - functional form)}
Let $f$ be a sufficiently nice function 
(for example, a function of the form
$f(x)=\int e^{ixy}\mu(dy)$ for a
complex measure $\mu$ with $\int\vert y\vert^2
\vert\mu\vert(dy)<\infty$). Then we have
$$df(S_t)=\partial f(S_t)\sharp dS_t+\frac 12
\Delta_t f(S_t)dt,$$
where $\partial f(X)$ is the extension of the derivation
$$\partial X^n=\sum_{k=0}^{n-1}X^k\otimes X^{n-k-1}$$
and $\Delta_t f$ denotes the function
$$\Delta_t f(x)=\frac\partial{\partial x}\int \frac{f(x)-f(y)}{x-y}
\nu_t(dy),$$
where $\nu_t$ is the distribution of $S_t$, i.e. a semicircular
distribution with variance $t$.

\begin{proof}
One has to check the statement for polynomials by using the 
product form of the Ito formula. All expressions in the statement
make also sense for nice functions and the statement extends by
continuity.
\\
For a polynomial $f(x)=x^n$ we have
$$dS_t^n=\sum_{k=0}^{n-1} S_t^kdS_tS_t^{n-k-1}+
\sum_{k,l\geq 0\atop k+l\leq n-2} S_t^kdS_tS_t^ldS_tS_t^{n-k-l-2}.$$
The first term gives directly the first term in our assertion
(this is just a non-commutative first
derivative), whereas the second yields
\begin{align*}
\sum_{k,l\geq 0\atop k+l\leq n-2} S_t^kdS_tS_t^ldS_tS_t^{n-k-l-2}
&=\sum_{k,l\geq 0\atop k+l\leq n-2} S_t^k\tau\lb S_t^l\rb dt
 S_t^{n-k-l-2}\\
&=
\sum_{k,l\geq 0\atop k+l\leq n-2} S_t^{n-l-2}\tau\lb S_t^l\rb dt\\
&=\sum_{m=1}^{n-1}m S_t^{m-1}\tau\lb S_t^{n-m-1}\rb dt,
\end{align*}
which can be identified with the second term in our statement.
One should note: the fact that the trace and not the identity acts on
the expression between two differentials
results finally in the unusual form of the second order
term in the functional form
of the Ito formula; 
it is not a non-commutative version of the second derivative,
but a mixture of derivative and difference expression.
\end{proof}

\bibliographystyle{amsplain}

\providecommand{\bysame}{\leavevmode\hbox to3em{\hrulefill}\thinspace}

\end{document}